\documentclass[a4paper]{amsart}
\usepackage[english]{babel}
\usepackage{amssymb}
\usepackage{graphicx}
\usepackage{booktabs}
\usepackage{longtable}
\usepackage[latin1]{inputenc}
\usepackage[T1]{fontenc}

\newcommand{\Z}{\mathbb{Z}}

\newcommand{\Q}{\mathbb{Q}}

\theoremstyle{definition}

\theoremstyle{remark}
\begin{document}
\title[Elliptic curves and torsion]{Elliptic curves with torsion group $\Z /6\Z $ }
\author[A. Dujella]{Andrej Dujella}
\address{Department of Mathematics\\University of Zagreb\\Bijeni{\v c}ka cesta 30, 10000 Zagreb, Croatia}
\email[A. Dujella]{duje@math.hr}

\author[J. C. Peral]{Juan Carlos Peral}
\address{Departamento de Matem\'aticas\\Universidad del Pa\'is Vasco\\Aptdo. 644, 48080 Bilbao, Spain}
\email[J. C. Peral]{juancarlos.peral@ehu.es}

\author[P. Tadi\' c]{Petra Tadi\' c}
\address{Department of Mathematics, Statistics and Information Science\\Juraj Dobrila University of Pula\\
Preradovi\' ceva 1, 52100 Pula, Croatia}
\email[P. Tadi\' c]{petra.tadic@unipu.hr}

\thanks{A. D. and P. T. were supported by the Croatian Science Foundation under the project no.
6422. J. C. P. was supported by the UPV/EHU grant EHU 10/05.}
\begin{abstract}

We exhibit several families  of elliptic curves with torsion group isomorphic to $ \Z/6\Z$ and generic rank at least $3$.
Families of this kind have been constructed previously by several authors:
Lecacheux, Kihara, Eroshkin and Woo.
We mention the details of some of them and we add other examples developed more recently by Dujella and Peral, and MacLeod.

Then we apply an algorithm  of Gusi\'c and Tadi\'c and we find the exact rank over $\Q(t)$ to be 3 and we also
determine free generators of the Mordell-Weil group for each family. By suitable specializations,
we obtain the known and new examples of curves over $\Q$ with torsion $ \Z/6\Z$ and rank $8$,
which is the current record.
\end{abstract}

\maketitle

\section{Introduction }
Families of curves with torsion group $\Z / 6 \Z$ and rank at least $3$ have been constructed previously by several authors, Lecacheux \cite{Le}, Kihara \cite{Ki} and later on Eroshkin \cite{Er} and Woo \cite{Wo}  are among them. More recently, Dujella and Peral (unpublished, 2012)
and MacLeod \cite{McL} found other examples.

We detail some of these families, and
in each case we include the relevant data and the output of the application of the algorithm of Gusi\'c and Tadi\'c \cite{GT}, to determine the exact rank over $\Q(t)$ and free generators of the Mordell-Weil group.
We apply \cite[Theorem 1.3]{GT}, which deals with elliptic curves $E$
given by $y^2 = x^3 + A(t)x^2 + B(t)x$, where $A,B \in \mathbb{Z}[t]$, with exactly
one nontrivial $2$-torsion point over $\mathbb{Q}(t)$.
If $t_0 \in \mathbb{Q}$ satisfies the condition that for every nonconstant
square-free divisor $h$ of $B(t)$ or $A(t)^2 - 4B(t)$ in $\mathbb{Z}[t]$
the rational number $h(t_0)$ is not a square in $\mathbb{Q}$,
then the specialized curve $E_{t_0}$ is elliptic and the specialization homomorphism at
$t_0$ is injective.
If additionally  there exist $P_1,\ldots,P_r\in E(\mathbb{Q}(t))$ such that
$P_1(t_0),\ldots,P_r(t_0)$ are  free generators of $E(t_0)(\mathbb{Q})$, then
$E(\mathbb{Q}(t))$ and $E(t_0)(\mathbb{Q})$ have the same rank $r$, and $P_1,\ldots,P_r$ are  free generators of $E(\mathbb{Q}(t))$.

Let us mention that the methods from \cite{GT} have been used in \cite{DGT}
to determine the exact rank over $\Q(t)$ of some of record rank families with torsion groups
$\mathbb{Z}/4\mathbb{Z}$,
$\mathbb{Z}/10\mathbb{Z}$, $\mathbb{Z}/12\mathbb{Z}$
and $\mathbb{Z}/2\mathbb{Z} \times \mathbb{Z}/8\mathbb{Z}$,
while the analogous result for torsion group
$\mathbb{Z}/2\mathbb{Z} \times \mathbb{Z}/4\mathbb{Z}$ has been obtained in \cite{DP1}
by using previous version of the injectivity criterium by Gusi\'c and Tadi\'c from \cite{GT1}.

\section{Curves with  torsion group $ \Z/6\Z$} \label{sec2}
The Tate's normal form for an elliptic curve is given by
\[
E(b,c): \quad y^2+(1-c) x y - b y = x^3- b x^2,
\]
see \cite[p. 147]{Kn}. It is nonsingular if and only if $b\neq 0$.
Using the addition law for  the point   $P=(0,0)$  we find
\begin{align*}
3 P&=(c, b-c),\\
- 3 P&=(c, c^2).
\end{align*}
It follows that $ P$ is a torsion point of   order $6$  for  $b=c+c^2$.
So, for $b=c+c^2$ we write the curve in the form $y^2= x^3+ A(c) x^2 + B(c) x.$ We have
\begin{align*}
A(c) &=1 + 6 c - 3 c^2\\
B(c)  &=-16 c^3
\end{align*}
In the new coordinates the torsion point of order $6$ is $(- 4 c, 4 c(1+c))$.

\section{Lecacheux' construction for rank $3$} \label{sec:le}

\subsection{Deduction of the coefficients for a rank $3$ family}\label{lecacheuxcurve}

A rank $3$ family comes from specialization in the rank $2$ family given by Lecacheux \cite{Le}.

In this paper Lecacheux constructs  four  biparametric families of curves with torsion group $\Z / 6 \Z$ and rank $2$ over $\Q(s,t)$, and by specialization of $s$ she gets  a family of curves with this torsion and rank 3 over $\mathbb Q(t)$. The starting model  used by Lecacheux for torsion  $\Z / 6 \Z$, once translated to  $(0, 0)$, is $y^2=x^3+ a_6(d) x^2+ b_6(d) x$ where
\begin{align*}
a_6(d)=& \, -3 - 6 d + d^2 , \\
b_6(d)=& \, 16  d.
\end{align*}

Lecacheux shows that using any of the four values
\begin{align*}
d_1=& \, \frac{(-t^2 + s) (1 + t s) (t - s^2)}{(-1 + t) t (-1 + s) s (t + s)},  \\
d_2=& \, \frac{(-s^2 + t (1 + t) (1 + s)) (-t^2 + (1 + t) s (1 + s))}{t s (1 + t +
   s) (t + s + t s)},\\
   d_3=& \, \frac{(t + s + t s (2 + t + s)) ((1 + t + s)^2 + t s (1 + t s + 2 (t + s)))}{t (1 + t + t^2) s (1 + s + s^2)},\\
   d_4=& \, -\frac{(t - 1) (s - 1) (t^2 + t s + s^2) (1 + t s + 2 (t + s))}{(t +
    s) (-t^2 + s) (t - s^2)},\\
\end{align*}
the resulting  biparametric family has rank $\geq 2$ over $\Q (s,t)$.  Then she specializes and gets a family of rank $3$ over $\Q (t)$. The value of $d$  is a rational expression  with numerator of degree $10$ and denominator of degree $9$ given by:
\[
\frac{(-1 + t + t^2) (-1 + 2 t + t^3) (-1 + 2 t - t^2 + t^3 + 2 t^4 + t^5)}{t^2 (1 + t) (-1 + 2 t) (-1 + 2 t + t^2) (-1 + t + t^2 + t^3)}.
\]

The curve of rank $3$ has polynomial coefficients in $t$ of degree $20$ and $40$ respectively. The coefficients are
{\footnotesize
\begin{eqnarray*}\label{ab63}
a_{63L}(t)&\!\!\!=\!\!\!&1 - 10 t + 35 t^2 - 36 t^3 - 61 t^4 + 146 t^5 + 5 t^6 - 254 t^7 +
 300 t^8 + 58 t^9 - 436 t^{10} - 294 t^{11} \\& & \mbox{} + 496 t^{12} + 710 t^{13 }+
 93 t^{14} - 434 t^{15} - 489 t^{16} - 264 t^{17} - 73 t^{18} - 6 t^{19} + t^{20} ,\\
b_{63L}(t)&\!\!\!=\!\!\!&16 t^6 (1 + t)^3 (-1 + 2 t)^3 (-1 + t + t^2) (-1 + 2 t + t^2)^3 (-1 +
   2 t + t^3) (-1 + t + t^2 + t^3)^3 \\ & & \mbox{}\times (-1 + 2 t - t^2 + t^3 + 2 t^4 +
   t^5).
\end{eqnarray*}
}
and the $x$-coordinates of three infinite order independent points  $P_1$, $P_2$ and $P_3$ are

{\footnotesize
\begin{align*}
x(P_1) = & \, 4 t (1 + t) (-1 + 2 t) (-1 + t + t^2) (-1 + 2 t + t^2)^2 (-1 + 2 t +
   t^3) (-1 + t + t^2 + t^3)^2,\\
x(P_2) = & \, 4 t^4 (1 + t) (-1 + 2 t) (-1 + 2 t + t^2)^2 (-1 + 2 t + t^3) (-1 +
   2 t - t^2 + t^3 + 2 t^4 + t^5),\\
x(P_3) = & \, 4 t^2 (-1 + 2 t) (-1 + 2 t + t^2)^2 (-1 + t + t^2 + t^3)^2 (-1 + 2 t -
    t^2 + t^3 + 2 t^4 + t^5).
\end{align*}}
\subsection{Generators}\label{firstproof}
Now we prove  that the elliptic curve $C$  over $\mathbb Q(t)$ given by the equation
$$C: \quad y^2=x^3+a_{63L}x^2+b_{63L}x, \ \ \mbox{(where $a_{63L},b_{63L}$ are given above)}$$
has  rank   equal  3 and the points $P_1,P_2, P_3$
are its free generators. We apply methods from \cite{GT}, as explained in
the introduction.

\begin{itemize}
\item We use the specialization  at $t_0=-7$ which is injective by \cite[Theorem 1.3]{GT}.
\item  The specialized curve $C^{-7}$ over $\mathbb Q$ with reduced coefficients is

$[0,298972834764046,0,-4129733728640949525768711375,0]$.

\item By {\tt mwrank} \cite{Cr}, the rank of this specialized curve over $\mathbb Q$  is
equal to 3 and its free generators are
\begin{eqnarray*}
G_1&\!\!\!=\!\!\!& [-242529733107900, -2078846881485761806650], \\
G_2&\!\!\!=\!\!\!& [-39216968008071, -749301352264289045760], \\
G_3&\!\!\!=\!\!\!& [31096987486425, 436763906196136192800].
\end{eqnarray*}
\item
We have that for the specialization of the points $P_1,P_2,P_3$ at $t_0=-7$ it holds
\begin{eqnarray*}
P_1^{-7}&\!\!\!=\!\!\!&   3\cdot T+0\cdot G_1 - 1\cdot G_2-1\cdot G_3,\\
P_2^{-7}&\!\!\!=\!\!\!&   1\cdot T+1\cdot G_1 + 0\cdot G_2-1\cdot G_3,\\
P_3^{-7}&\!\!\!=\!\!\!&   5\cdot T+0\cdot G_1 -1\cdot G_2+0\cdot G_3,
\end{eqnarray*}
where $T$ is the torsion point of order 6 on the specialized curve.
\item
Since the determinant of the base change matrix has  absolute value 1, it follows that the points $P_1, P_2 , P_3$  with $x$-coordinates  given above are free generators of  $C$ over $\mathbb Q(t)$.
\end{itemize}

% \subsection{Examples of high rank}

\section{Kihara's construction for rank $3$} \label{sec:ki}

\subsection{Family data}\label{kiharacurve}
Kihara, \cite{Ki}, consideres the projective curve
\[
 (x^2-y^2)^2+ ( 2 a x^2+ 2 b y^2) z^2+ c z^4=0,
\]
then he uses the  substitutions
\begin{align*}
X=& \, \frac{x^2}{y^2},\\
Y=& \, \frac{x(c z^2+ a x^2+ b y^2)}{y^3},
\end{align*}
and arrives to  the elliptic curve
\[
E:\quad Y^2=X((a^2- c) X^2+ (2 a b + 2 c) X+ (b^2-c)).
\]
The point $P=(1, a+ b)$ is on $E$. Forcing $P$ to be of order $3$ implies $c=\frac{(3a -b)(a+b)}{4}$ in which case $P+(0,0)$ is of order $6$. With this value of $c$ the curve $E$ is given by $Y^2= X^3+ A(a,b)  X^2  + B(a,b) X$ where
\begin{align*}
A(a,b)=&  \, 2 (3 a^2 + 6 a b - b^2),   \\
B(a,b)=& \, -(a - b)^3 (3 a + 5 b).
\end{align*}
By imposing new points on the curve,  Kihara gets a family with torsion group $\Z/6\Z$ and rank at least 3 given as $y^2=x^3+ a_{63K}(t) x^2+ b_{63K}(t) x$ where
\begin{align*}
a_{63K}(t)=& \, -2(64 t^8 - 1952 t^7 - 4652 t^6 - 10172 t^5 - 28955 t^4 +
   35602 t^3 \\& \mbox{}- 56987 t^2 + 83692 t + 9604)      \\
b_{63K}(t)=&\, (t - 7)^3 (t + 2)^3 (2 t + 1)^3 (4 t - 7)^3 (2 t^2 - 91 t +
   98) (4 t^2 + 13 t + 1)
\end{align*}

The $x$-coordinates of three  independent points of infinite order are:
\begin{align*}
x(Q_1)=&\frac{9 (-7 + t)^2 (2 + t)^2 (1 + 2 t)^2 (-7 + 4 t)^2 (7 + 2 t^2)^2}{(-7 - 4 t + 2 t^2)^2},\\
x(Q_2)=&\frac{(-7 + t)^2 (2 + t)^2 (1 + 2 t)^2 (-7 + 4 t)^2 (7 + 2 t^2)^2}{(-7 + 14 t + 2 t^2)^2},\\
x(Q_3)=&\frac{(-7 + t)^2 (2 + t)^2 (1 + 2 t)^2 (-7 + 4 t)^2 (-7 + 14 t + 2 t^2)^2}{(-7 - 4 t + 2 t^2)^2}.
\end{align*}

\subsection{Generators}
Now following the modified proof of subsection \ref{firstproof} we prove that Kihara's elliptic curve  over $\mathbb Q(t)$ given at the end of subsection \ref{kiharacurve} has rank  equal to 3 and the points $W_1,W_2,Q_3$ are its free generators  (the $x$-coordinate of $Q_3$ given above, and of $W_1, W_2$ given below). So the above given points $Q_1,Q_2,Q_3$ found by Kihara, are not its free generators.
\begin{itemize}
\item The specialization  at $t_0=15$ is injective by \cite[Theorem 1.3]{GT}.

\item The reduced coefficients  of the specialized curve are

$[0,7236353038,0,-2438945400771712139,0].$

\item By {\tt mwrank} \cite{Cr}, the specialized curve has rank equal to $3$ and its free generators are
\begin{eqnarray*}
G_1&\!\!\!=\!\!\!&[-3875031,-3091858818006],\\
G_2&\!\!\!=\!\!\!&[5184423886489/64, -12318461250340505085/512],\\
G_3&\!\!\!=\!\!\!&[83648786841/4, -28010191585888407/8].
\end{eqnarray*}
\item We have that the specialization of the points $Q_1,Q_2,Q_3$ at $t_0=15$ satisfies
\begin{eqnarray*}
Q_1^{15}=0\cdot T+  2\cdot G_1+  0\cdot G_2+  1\cdot G_3,\\
Q_2^{15}=0\cdot T+  2\cdot G_1  +2\cdot G_2 + 3\cdot G_3,\\
Q_3^{15}=4\cdot T + 0\cdot G_1  +0\cdot G_2  -1\cdot G_3,
\end{eqnarray*}
here $T$ is the torsion point of order 6 on the specialized curve.
\item From this we see that we can take points
   $W_1$ and $W_2$  such that
\begin{eqnarray*}Q_1+Q_3&=&2W_1,\\
Q_2+3Q_3&=&2W_2.
\end{eqnarray*}
\item The
$x$-coordinates of $W_1$  and $W_2$ are
\begin{eqnarray*}
x(W_1)&\!\!\!=\!\!\!&(4t-7)(t-7)(4t^2+13t+1)(2t^2-7t+14)^2,\\
x(W_1)&\!\!\!=\!\!\!&(2t+1)^2(t+2)^2(4t^2+13t+1)(64t^5-536t^4+1324t^3+224t^2-490t-343)^2\\
& & \mbox{}\times (t-7)(4t-7)/(64t^5-8t^4-284t^3+44t^2+476t-49)^2,
\end{eqnarray*}

\item Now it is obvious that   $W_1, W_ 2,  Q_3$ are free generators   of the observed elliptic curve over $\mathbb Q(t)$. More precisely
\begin{eqnarray*}
W_1^{15}=1\cdot T-1\cdot G_1+  0\cdot G_2+  0\cdot G_3\\
W_2^{15}=3\cdot T+  1\cdot G_1  +1\cdot G_2 + 3\cdot G_3\\
Q_3^{15}=4\cdot T + 0\cdot G_1  +0\cdot G_2  -1\cdot G_3.
\end{eqnarray*}

\end{itemize}

%\subsection{Examples of high rank}

\section{Eroshkin's  construction for rank $3$} \label{sec:er}

\subsection{Family data}\label{eroshkincurve}

In 2008, Eroshkin \cite{Er} constructed another example of rank $3$ curve over $\Q (t)$.
The starting point for his construction is the two-parametric family with rank $\geq 2$
over $\Q (u,v)$, given by
\begin{eqnarray*}
A(u,v) &\!\!\!=\!\!\!&
v^8u^4+2u^7v^5+3u^6v^6+2u^5v^7+3u^4v^2+3v^4u^2-12u^3v^3+u^8v^4-3u^6 \\
& & \hspace{-40pt}\mbox{}-6u^5+7u^4-3u^2+6u^3-6v^5+7v^4-3v^6+21u^6v^2+36u^4v^4-9u^2v^2 \\
& & \hspace{-40pt}\mbox{}+12v^2u-24u^4v^5-9u^4v^6+36u^5v^3+21u^2v^6-24u^5v^4-10v^4u-6v^7u^2 \\
& & \hspace{-40pt}\mbox{}+4u^5v^6-2u^2v^5+2u^4v^7-14u^3v^6+8u^5v-9u^6v^4-18u^5v^5+2v^7u^3+12vu^2 \\
& & \hspace{-40pt}\mbox{}+2u^3v-6v^6u-20u^2v^3+8v^5u-6u^6v+36u^3v^5+2u^7v^3-10u^4v+4u^6v^5 \\
& & \hspace{-40pt}\mbox{}-6uv+2v^3u-6u^7v^2-14u^6v^3-3v^2+6v^3+2u^7v^4-2u^5v^2-20u^3v^2,   \\
B(u,v) &\!\!\!=\!\!\!&
-16(u+v)^3(-1+v)^3(1+v)^3(-1+u)^3(1+u)^3 \\
& & \hspace{-40pt}\mbox{}\times (u^2+uv+v^2-2vu^2-2v^2u+u^3v+v^3u-2u^3v^2-2u^2v^3+u^4v^2+v^4u^2+u^3v^3).
\end{eqnarray*}
The curve $y^2=x^3+A(u,v)x^2+B(u,v)x$ contains the points with $x$-coordinates
$4(v-1)^2(1+v)^3(-1+u)^2(1+u)^3$ and $-4(u-1)(u+1)(v-1)^3(v+1)^2(u+v)^3$.
Forcing the point with $x$-coordinate $3(v-1)^2(v+1)^2(-1+u)^2(1+u)^2(v+u)^2$ to belong to the curve,
leads to the condition $v=1/3$. By substitution $u=(1-t)/(1+t)$ and some simplifications,
we get Eroshkin's rank $3$ curve
with coefficients
\begin{align*}
a_{63E}(t)=& \, 16 + 576 t - 1408 t^2 - 1440 t^3 + 1608 t^4 + 720 t^5 - 352 t^6 -
 72 t^7 + t^8,    \\
b_{63E}(t)=& \, 27648 (-2 + t)^3 t^3 (1 + t)^3 (2 + t^2)^2,
\end{align*}
and $x$-coordinates of three independent points $R_1,R_2,R_3$ given by:
\begin{align*}
x(R_1)=& \, 864 (-2 + t)^3 t^3,\\
x(R_2)=& \, 3456 t^2 (1 + t)^3,\\
x(R_3)=& \, 288 (-2 + t)^3 t^3 (1 + t)^2.
\end{align*}

Another way to get a rank $3$ family from Eroshkin's two-parametric family is to take 
$v=-1/3$. By substitution $u=(1-t)/(1+t)$, as above, we get 
\begin{align*}
a_{63ER}(t)=& \, 256t^8-2304t^7-3232t^6+1008t^5+2337t^4-504t^3-808t^2+288t+16,\\
b_{63ER}(t)=& \, 27648(16t^4-11t^2+4)(2t-1)^3(t+1)^3t^3,
\end{align*}
and three infinite order independent points  $S_1,S_2,S_3$ with $x$-coordinates 
\begin{align*}
x(S_1) = & \,1728t^2(t+1)^3,\\
x(S_2) = & \, 864t^3(2t-1)^3,\\
x(S_3) = & \,864t^3(t+1)^2(2t-1).
\end{align*}

%%%%%%%%%%%%%%%%%%%%%%%%

\subsection{Generators}
We prove that the rank of both families (for $v=1/3$ and $v=-1/3$) is equal to 3 and the points 
$R_1,R_2, R_3$ (resp. $S_1,S_2,S_3$) are free generators of the elliptic curve 
$y^2=x^3+a_{63E}x^2+b_{63E}x$ (resp. $y^2=x^3+a_{63ER}x^2+b_{63ER}x$) over $\mathbb Q(t)$, 
where $a_{63E},b_{63E},a_{63ER},b_{63ER}$ are given above.
\begin{itemize}
\item For the first curve, we use the specialization at  $t_0=-11$ 
which is injective by \cite[Theorem 1.3]{GT}.
\item The specialized curve over $\mathbb Q$ with reduced coefficients is
$$ [0,100352409,0,-15100698522624000,0]. $$
\item 
By {\tt mwrank}\cite{Cr}, the rank of this specialized curve over $\mathbb Q$  is
equal to 3 and its free generators are
\begin{eqnarray*}
G_1&\!\!\!=\!\!\!&[-11073920,-422104608640] \\
G_2&\!\!\!=\!\!\!&[-21632000,-602905472000], \\
G_3&\!\!\!=\!\!\!&[879994086400/9409, -486030413853670400/912673];
\end{eqnarray*}
\item
For the specialization of the points $R_1,R_2,R_3$ at $t_0=-11$ we have
\begin{eqnarray*}
R_1^{-11}&\!\!\!=\!\!\!&   5\cdot T+0\cdot G_1 + 1\cdot G_2+1\cdot G_3,\\
R_2^{-11}&\!\!\!=\!\!\!&   5\cdot T+1\cdot G_1 -1\cdot G_2-1\cdot G_3,\\
R_3^{-11}&\!\!\!=\!\!\!&   0\cdot T-1\cdot G_1 +1\cdot G_2+0\cdot G_3,
\end{eqnarray*}
where $T$ is the torsion point of order 6 on the specialized curve.
\item
Since the determinant of the base change matrix has  absolute value 1, it follows that the points $R_1, R_2 , R_3$  with $x$-coordinates  given above are free generators of  
$y^2=x^3+a_{63E}x^2+b_{63E}x$ over $\mathbb Q(t)$.
\end{itemize}

\medskip

\begin{itemize}
\item For the second curve we use the specialization at $t_0=-15$ 
which is injective by \cite[Theorem 1.3]{GT}.
\item  The specialized curve  over $\mathbb Q$ is
$$[0,1012299875521,0,-6159774508321416192000,0].$$
\item By {\tt mwrank} \cite{Cr}, the rank of this specialized curve over $\mathbb Q$  is
equal to 3 and its free generators are
{\footnotesize{
\begin{eqnarray*}
G_1&\!\!\!=\!\!\!&[4710284355840,11266618646136456960], \\
G_2&\!\!\!=\!\!\!&[-70907506432,-71902551469760256], \\
G_3&\!\!\!=\!\!\!&[\frac {21570079809600}{169},-\frac{292982182574714539200}{2197}].
\end{eqnarray*}
}}
\item
We have that for the specialization of the points $S_1,S_2,S_3$ at $t_0=-15$ it holds
\begin{eqnarray*}
S_1^{-15}&\!\!\!=\!\!\!&   1\cdot T-1\cdot G_1 +0\cdot G_2-1\cdot G_3,\\
S_2^{-15}&\!\!\!=\!\!\!&   3\cdot T+0\cdot G_1 +1\cdot G_2+0\cdot G_3,\\
S_3^{-15}&\!\!\!=\!\!\!&   1\cdot T+0\cdot G_1 +1\cdot G_2-1\cdot G_3,
\end{eqnarray*}
where $T$ is the torsion point of order 6 on the specialized curve.
\item
Since the determinant of the base change matrix has  absolute value 1, it follows that the points $S_1, S_2 , S_3$  with $x$-coordinates  given above are free generators of 
$y^2=x^3+a_{63ER}x^2+b_{63ER}x$ over $\mathbb Q(t)$.
\end{itemize}

%\subsection{Examples of high rank}

%%%%%%%%%%%%%%%%%%%%%%%%%%%%%%%%%%%%%%
\section{A direct construction by Dujella and Peral (2012)} \label{sec:dp}

\subsection{Deduction of the family}

Our next example is based on the use of $2$-descent as follows, see e.g. \cite{DP}.
Consider the general curve with  torsion $\Z / 6 \Z$  as given in Section \ref{sec2}
i.e.: $$y^2= x^3+ x^2  (1 + 6 c - 3 c^2) + x (-16 c^3).$$

We first impose $16 \, c^2$ as $x$-coordinate of a new point, this is equivalent  to solve $ 1+5 c+13 c^2=\square$. The corresponding parametrization is  $c = -\frac{(-4 + u) (-2 + u)}{-13 + u^2}$.  This gives us an infinite order point and a curve  with rank at least $1$ over $\Q (u)$ whose coefficients are
\begin{align*}
A_{61}(u)=& \, 601 - 180 u - 152 u^2 + 72 u^3 - 8 u^4 , \\
B_{61}(u)=& \, 16 (-4 + u)^3 (-2 + u)^3 (-13 + u^2).
\end{align*}

Now observe that imposing  $4 (-4 + u) (-2 + u)^3$ as a new point is the same as  to   solve $-103 + 12 u + 4 u^2 =\square$. This can be achieved  with $u=-\frac{103 + t^2}{4 (-3 + t)}$.
The new family has the following  coefficients
\begin{align*}
A_{63}(t)=& \, -2 (26009437 + 18059772 t + 5057576 t^2 + 813612 t^3 + 89370 t^4\\& +
   7860 t^5 + 608 t^6 + 36 t^7 + t^8), \\
B_{63}(t)=& \, (5 + t)^3 (11 + t)^3 (79 + 8 t + t^2)^3 (8737 + 1248 t - 2 t^2 + t^4).
\end{align*}
In this case the  family has rank $\geq 3$.
Below we list the  $x$-coordinates of three independent points of infinite order on the curve
\begin{align*}
x(P_1)=& \, 4 (5 + t)^2 (11 + t)^2 (79 + 8 t + t^2)^2,\\
x(P_2)=& \, (5 + t) (11 + t) (79 + 8 t + t^2)^3,\\
x(P_3)=& \, \frac{64 (4 + t)^2 (5 + t)^2 (11 + t)^2 (79 + 8 t + t^2)^2}{(29 + 6 t +
  t^2)^2}.
\end{align*}
The first two points are the ones we have  imposed in our constructions, while the  third one  appears in the last change. So this change produces two new independent points of infinite order and accordingly this family has rank at least $3$.

\subsection{Generators}

We prove the curve  $y^2= x^3+ A_{63}x^2+ B_{63}x$ has rank exactly $3$ over $\Q (t)$ with free generators the points $P_1,P_2,P_3$.

Here is the proof, it is similar to the proof in subsection \ref{firstproof}.
\begin{itemize}
\item We take the  specialization homomorphism at $t_0 = 13$ which is injective by \cite[Theorem 1.3]{GT}.
%\item The
%transformation for the reduced form is $ 1/96$.
\item  Reduced form of the specialized curve is  $[0, -3121367, 0, 2201786966016, 0]$.
\item
{\tt Mwrank} \cite{Cr} shows that the rank of the specialized curve is equal to $3$, and its free generators are
\begin{eqnarray*}
G_1&\!\!\!=\!\!\!& [38016, 281508480],\\
G_2&\!\!\!=\!\!\!& [213276816, -3091894330320],\\
G_3&\!\!\!=\!\!\!& [381655296/361, 985650589440/6859].
\end{eqnarray*}
We have that for the specialization of the points $P_1,P_2,P_3$ at $t_0=13$
\begin{eqnarray*}
P_1^{-11}&\!\!\!=\!\!\!&   0\cdot T-2\cdot G_1 + 1\cdot G_2+0\cdot G_3,\\
P_2^{-11}&\!\!\!=\!\!\!&   0\cdot T-1\cdot G_1+ 0\cdot G_2-1\cdot G_3,\\
P_3^{-11}&\!\!\!=\!\!\!&   2\cdot T+2\cdot G_1 -1\cdot G_2+1\cdot G_3,
\end{eqnarray*}
where $T$ is the torsion point of order 6 on the specialized curve.
\item
The points with $x$-coordinates $P_1(t), P_2(t), P_3(t)$ are free generators of the full Mordell-Weil group, analogously as before.
\end{itemize}
%\subsection{Examples of high rank}

\section{The construction by MacLeod (2013)} \label{sec:ml}
\subsection{Coefficients for a rank $3$ family}\label{macleodcurve}

A rank $3$ family comes from the paper of MacLeod.
The curve of rank $3$ has polynomial coefficients:
{\footnotesize
\begin{eqnarray*}\label{abML}
a_{ML}(t)&=&60637t^{16}-9446576t^{15}+318481560t^{14}+7559095920t^{13}-673080292884t^{12}\\
&&+13070347480656t^{11}-22695122540632t^{10}-1545816592882960t^9\\
&&-5109785850000978t^8+496074059556450416t^7-4896348487968033496t^6\\
&&+8785946885865627600t^5+127778653981127482476t^4-861560195308301691408t^3\\
&&+2072099226498542082072t^2-1499601599678467324208t-1357666940926062868067,\\
b_{ML}(t)&=&-(t^2-14t-83)^3(17t^2-226t-367)^3(t^2+58t-347)^3(t^2-290t+2017)^3\times\\
&&\times(445t^8-26680t^7+298540t^6+7487800t^5-100559026t^4-816359048t^3\\
&&+14922867436t^2-55316709112t+68821828189),\end{eqnarray*}
}
and the $x$-coordinates of three infinite order independent points  $Q_1$, $Q_2$ and $Q_3$ are
{\footnotesize
\begin{align*}
x(Q_1) = & \,-4(13t^2-458t+1021)^2(t^2+4t-149)^2(t^2-290t+2017)(t^2+58t-347)\times\\
&\times(17t^2-226t-367)(t^2-14t-83),\\
x(Q_2) = & \, -9(17t^2-226t-367)^3(t^2-14t-83)^3(t^2-290t+2017)(t^2+58t-347),\\
x(Q_3) = & \,(445t^8-26680t^7+298540t^6+7487800t^5-100559026t^4-816359048t^3+14922867436t^2\\
&-55316709112t+68821828189)(17t^2-226t-367)^2(t^2-14t-83)^2.
\end{align*}
}
\subsection{Generators}
Following the proof in subsection \ref{firstproof} we prove that MacLeod's elliptic curve over $\mathbb Q(t)$ given in subsection \ref{macleodcurve} has rank   equal to 3 over $\Q(t)$ and the points $Q_1,Q_2, Q_3$
are its free generators. 
\begin{itemize}
\item We use the specialization  at $t_0=-30$ which is injective by \cite[Theorem 1.3]{GT}.
\item  The specialized curve  over $\mathbb Q$  already has the reduced coefficients  
{\footnotesize{
$$[0,-396608000022758283110679701027,0,$$
$$440417648413441630018736813558974274749049
48473305010381651,0].$$
}}
\item By {\tt mwrank} \cite{Cr}, the rank of this specialized curve over $\mathbb Q$  is
equal to 3 and its free generators are
{\footnotesize{
\begin{eqnarray*}
G_1&\!\!\!=\!\!\!&[1371227854300350505685906632351,1375839998833531665451365906419334737742953475], \\
G_2&\!\!\!=\!\!\!&\Big[\frac{11719103414610661317429257961549638649}{259596544},\\
&\!\!\!\!\!\!&-\frac{149169932197111418433971723731890785699819555746182187635}{4182619516928}\Big], \\
G_3&\!\!\!=\!\!\!&\Big[\frac {3986901670582496234640210004431381477056020081}{33036918984050881},\\
&\!\!\!\!\!\!&-\frac {216209573005548181556962517949975521622133551319748523058488502950375}{6004810514024749417546529}\Big].
\end{eqnarray*}
}}
\item
We have that for the specialization of the points $Q_1,Q_2,Q_3$ at $t_0=-30$ it holds
\begin{eqnarray*}
Q_1^{-30}&\!\!\!=\!\!\!&   0\cdot T+1\cdot G_1 +1\cdot G_2-1\cdot G_3,\\
Q_2^{-30}&\!\!\!=\!\!\!&   2\cdot T+1\cdot G_1 + 0\cdot G_2-1\cdot G_3,\\
Q_3^{-30}&\!\!\!=\!\!\!&   1\cdot T-1\cdot G_1 +0\cdot G_2+0\cdot G_3,
\end{eqnarray*}
where $T$ is the torsion point of order 6 on the specialized curve.
\item
Since the determinant of the base change matrix has  absolute value 1, it follows that the points $Q_1, Q_2 , Q_3$  with $x$-coordinates  given above are free generators of the observed curve over $\mathbb Q(t)$.
\end{itemize}

\section{Examples with rank 8 over $\Q$}  \label{sec:ex}

Currently the highest know rank for curves over $\Q$ with torsion group $\Z /6\Z $ is $8$.
In this section, we explain how 10 previously known and 10 new curves with rank $8$ were (re)discovered.
The details (minimal equation, torsion points, independent points of infinite order)
for all these curve can be found at $\Z /6\Z $ section of tables \cite{Du}.

We applied standard techniques for finding curves with high rank in families
with reasonably large generic rank. These techniques include computations of
Mestre-Nagao sums \cite{Me,Na} and 2-Selmer rank. In our computations we use PARI/GP \cite{P}
and {\tt mwrank} \cite{Cr}. In most of the families mentioned in previous sections,
we are able to find curves with rank at least $7$. Here we report only constructions
which give curves with rank $8$.

\subsection{Kihara's family}

By applying the above mentioned techniques to Kihara's family of rank 3 over $\Q(t)$
from Section \ref{sec:ki}, we find the curves with rank 8 over $\Q$ for the following
values of the parameter $t$:
\begin{itemize}
\item $39/121$ (2-isogenous to 4th curve in \cite{Du}; Elkies (2008)),
\item $178/113$ (2-isogenous to 5th curve in \cite{Du}; Elkies (2008)),
\item $527/301$ (8th curve in \cite{Du}; Dujella (2008)),
\item $-1/839$ (9th curve in \cite{Du}; Dujella (2008)),
\item $519/407$ (7th curve in \cite{Du}; Dujella (2008)),
\item $722/211$ (6th curve in \cite{Du}; Elkies (2008)),
\item $-601/383$ (11th curve in \cite{Du}; Dujella-Peral-Tadi\'c (2014)),
\item $-909/652$ (12th curve in \cite{Du}; Dujella-Peral-Tadi\'c (2014)).
\end{itemize}

\subsection{Family by Dujella and Peral}
In the family of rank 3 by Dujella and Peral from Section \ref{sec:dp}
we found only one curve with rank $8$, for $t=-244/3$. It is 2-isogenous
to 10th curve in \cite{Du}, listed there as Dujella-Peral (2012).

\subsection{Eroshkin's two-parametric family}

In Eroshkin's family of rank 3 over $\Q(t)$ obtained for $v=1/3$ we were able to find 
only examples with rank 7. 
On the other hand, in the family of rank $3$ obtained for $v=-1/3$ we have found four examples 
with rank 8. 
Moreover, several examples with rank 8 were found within 
Eroshkin's two-parametric family with rank $\geq 2$ over $Q(u,v)$. 
Here we list parameters $(u,v)$ corresponding to these examples:
\begin{itemize}
\item $(20/19,  -7/17)$ (1st curve in \cite{Du}; Eroshkin (2008)),
\item $(-4/25,  7/19)$ (2nd curve in \cite{Du}; Eroshkin-Dujella (2008)),
\item $(-73/80,  11/7)$ (3rd curve in \cite{Du}; Eroshkin-Dujella (2008)),
\item $(-1/31,  21/13)$ (13th curve in \cite{Du}; Dujella-Peral-Tadi\'c (2014)),
\item $(101/119,  -20/17)$ (14th curve in \cite{Du}; Dujella-Peral-Tadi\'c (2014)),
\item $(77/338,  -1/3)$ (15th curve in \cite{Du}; Dujella-Peral-Tadi\'c (2015)),
\item $(633/229,  -1/3)$ (16th curve in \cite{Du}; Dujella-Peral-Tadi\'c (2015)),
\item $(-283/437,  6/19)$ (17th curve in \cite{Du}; Dujella-Peral-Tadi\'c (2015)),
\item $(-299/1395,  -1/3)$ (18th curve in \cite{Du}; Dujella-Peral-Tadi\'c (2015)),
\item $(109/466,  -7/8)$ (19th curve in \cite{Du}; Dujella-Peral-Tadi\'c (2015)),
\item $(447/305,  -1/3)$ (20th curve in \cite{Du}; Dujella-Peral-Tadi\'c (2015)).
\end{itemize}

\end{document}